\documentclass{jloganal}

\usepackage{pinlabel}
\usepackage{amssymb}
\usepackage{amscd}
\usepackage{amsfonts}
\usepackage{mathrsfs}
\usepackage{setspace}
\usepackage{mathtools}

\title[Transfer principle]{A transfer principle for second-order arithmetic, \\and applications}
\author[M Carl]{Merlin Carl}
\givenname{Merlin}
\surname{Carl}
\address{Department of Mathematics and Statistics\\
University of Konstanz \\\newline
78457 Konstanz\\DE}
\email{merlin.carl@uni-konstanz.de}
\urladdr{}
\author[A Jamneshan]{Asgar Jamneshan}
\givenname{Asgar}
\surname{Jamneshan}
\address{Department of Mathematics\\ 
ETH Zurich \\\newline 
R\"amistrasse 101\\ 8092 Zurich \\ CH\\\newline 
A.J.~gratefully acknowledges financial support from DFG project KU 2740/2-1.} 
\email{asgar.jamneshan@math.ethz.ch}
\urladdr{}

\arxivreference{arXiv:1804.01573}
\arxivpassword{4us27}

\numberwithin{equation}{section} 
\theoremstyle{plain}
\newtheorem{theorem}[subsection]{Theorem}
\newtheorem{proposition}[subsection]{Proposition}
\newtheorem{lemma}[subsection]{Lemma}
\newtheorem{corollary}[subsection]{Corollary}
\newtheorem{definition}[subsection]{Definition}
\newtheorem{remark}[subsection]{Remark}

\newtheorem{examples}[subsection]{Examples}

\renewcommand{\leq}{\leqslant}
\renewcommand{\geq}{\geqslant}

\renewcommand\R{\mathbb{R}}
\newcommand\N{\mathbb{N}}
\newcommand\F{\mathcal{F}}

\begin{document}

\begin{abstract}
In the theory of conditional sets, many classical theorems from areas such as functional analysis, probability theory 
or measure theory are lifted to a conditional framework, often to be applied in areas such as 
mathematical economics or optimization. The frequent experience that
such theorems can be proved by `conditionalizations' of the classical proofs suggests that a general 
transfer principle is in the background, and that formulating and proving such a transfer principle would 
yield a wealth of useful further conditional versions of
classical results, in addition to providing a uniform approach to the results already known. In this paper, 
we formulate and prove such a transfer principle based on second-order arithmetic, which, by the results of 
reverse mathematics, suffices for the bulk of
classical mathematics, including real analysis, measure theory and countable algebra, and excluding only more 
remote realms like category theory, set-theoretical topology or uncountable set theory, see e.g. the introduction of Simpson \cite{simpson2009subsystems}.
This transfer principle is then employed to give short and easy proofs of conditional versions of central results 
in various areas of mathematics, including theorems which have not been proven by hand previously such as Peano existence theorem, Urysohn's lemma and the Markov-Kakutani fixed point theorem. 
Moreover, we compare the interpretation of certain structures in a conditional model with their meaning in a standard model. 
\end{abstract}

\begin{asciiabstract}
In the theory of conditional sets, many classical theorems from areas such as functional analysis, probability theory 
or measure theory are lifted to a conditional framework, often to be applied in areas such as 
mathematical economics or optimization. The frequent experience that
such theorems can be proved by `conditionalizations' of the classical proofs suggests that a general 
transfer principle is in the background, and that formulating and proving such a transfer principle would 
yield a wealth of useful further conditional versions of
classical results, in addition to providing a uniform approach to the results already known. In this paper, 
we formulate and prove such a transfer principle based on second-order arithmetic, which, by the results of 
reverse mathematics, suffices for the bulk of
classical mathematics, including real analysis, measure theory and countable algebra, and excluding only more 
remote realms like category theory, set-theoretical topology or uncountable set theory, see e.g. the introduction of Simpson \cite{simpson2009subsystems}.
This transfer principle is then employed to give short and easy proofs of conditional versions of central results 
in various areas of mathematics, including some theorems which have not been proven by hand previously such as Peano existence theorem, Urysohn's lemma and the Markov-Kakutani fixed point theorem. 
Moreover, we compare the interpretation of certain structures in a conditional model with their meaning in a standard model. 
\end{asciiabstract}

\maketitle

\section{Introduction}

Fixing a probability space $(\Omega,\mathcal{F},\mathbb{P})$, one can distinguish between probabilistic and deterministic objects such as a deterministic real number which is an element of $\mathbb{R}$ and a random real number which is a measurable function $x\colon \Omega \to \mathbb{R}$. 
Extending such reasoning, one might speak of a `random', `stochastic', 'measurable' or `conditional' version of a `deterministic', `classical' or `standard' theorem which expresses a randomization of its statement.  
 To illustrate, a conditional version of the Bolzano-Weierstra\ss~theorem states that for every sequence $(x_k)$ of random real numbers such that $x=\limsup x_k<\infty$ almost surely, there exists a strictly increasing sequence $n_1<n_2<\ldots$ of integer-valued random variables such that $x_{n_k}$ converges almost surely to $x$.  
Experience has shown that many classical theorems have such conditional analogues such as the Heine-Borel theorem, the Hahn-Banach extension and separation theorems and the Brouwer fixed point theorem, see e.g.~Cheridito et al.~\cite{cheridito2015conditional}, Drapeau et al.~\cite{drapeau2016algebra,drapeau2013brouwer}, Filipovic et al.~\cite{filipovic2009separation}, and Jamneshan et al.~\cite{jamneshan2017measures,jamneshan2017compact} for an account. 
It is thus tempting to aim for a general transfer principle that allows one to `import' classical theorems into a conditional setting. 

In section \ref{s3}, we prove such a transfer principle
based on second-order logic. Strong arguments have been put forward in favor of the claim that second-order logic is a satisfying formal framework for the bulk
of classical mathematics, see e.g.~the introduction to Simpson \cite{simpson2009subsystems}; and this was impressively confirmed by the results of reverse mathematics. 
More precisely, we prove that any consequence of the second-order axiomatic system of arithmetical comprehension ACA$_{0}$, which has a second-order comprehension axiom for formulas in which all quantifiers range
over natural numbers (see Simpson \cite{simpson2009subsystems}), also holds conditionally. 
To this end, we show (i) that the axioms of ACA$_{0}$ hold in the structure with first order part $L^{0}(\mathbb{N})$ and second-order part its conditional power set $\mathcal{P}$ with truth value $\Omega$, and 
(ii) that truth value $\Omega$ is preserved by the usual deduction rules of second-order predicate calculus. 
To this end, a conditional element relation between $L^0(\N)$ and $\mathcal{P}$ is introduced.  
We expect that, with a certain amount of extra technical effort,
this can be extended to a considerably stronger transfer principle for full second-order logic.

In section \ref{s4}, we then discuss consequences of the transfer principle. 
We verify that the transfer principle yields a conditional version of a whole class of classical theorems  many of which were proved by hand previously\footnote{The latter practice  provides  some useful insight into the transfer process.}. 
As new consequences, we obtain a conditional version of the Peano existence theorem, Urysohn's lemma, the existence of an orthonormal basis and the Markov-Kakutani theorem.

A conditional version of a classical theorem is oftentimes also a theorem about a more involved situation in a classical setting.  
For instance, the above conditional version of the Bolzano-Weierstra{\ss}-theorem is also a classical theorem securing existence of measurably parametrized almost surely converging subsequences of an almost surely upper bounded sequence of real-valued measurable functions\footnote{ 
In this classical form under the name of a `randomized', `stochastic' or `measurable' version of the Bolzano-Weierstra{\ss}-theorem 
this statement is proved in F\"ollmer and Schied \cite{foellmer2011stochastic} and Kabanov and Stricker \cite{kabanov2001teachers}, motivated by applications in mathematical economics.}.  
We will systematically investigate a standard interpretation of a conditional version of a classical structure or theorem in section \ref{s4}. 
As a result, bridges to the analysis in $L^0(\R)$-modules (e.g.~Cheredito et al.~\cite{cheridito2015conditional}, Filipovic et al.~\cite{filipovic2009separation}) and set-valued analysis (e.g.~Rockafellar and Wets \cite{rockafellar02}) are built. 

We will be also interested in the reverse direction. 
Namely, it can be verified that certain classical structures have a useful interpretation in a conditional model.  
For example, the standard space $L^0(\mathbb{R})$ of real-valued Borel measurable functions on a $\sigma$-finite measure space $(\Omega,\mathcal{F},\mu)$ modulo almost everywhere equality are the real numbers in a conditional 
model.  
More generally, if $E$ is a complete separable metric space and $L^0(E)$ is the space of $E$-valued Borel measurable functions on $\Omega$ modulo almost everywhere equality, 
we show that $L^0(E)$ can be identified with a complete separable metric space in a conditional model of ACA$_{0}$.  
The space $L^0(E)$ reflects a measurable parametrization of the elements of $E$ relative to a base space $(\Omega,\mathcal{F},\mu)$.  
Such a measurable parametrization is a constituent part of a conditional model\footnote{Therefore conditional models might equally be called `stochastic' or `measurable' models of an axiomatic system.}. 
From an external perspective, one may view the transfer principle as a device which parametrizes classical theorems in a measurable way relative to a fixed measure space. 
In particular, the application of a conditional version of classical theorems preserves measurability, and thus provides an 
alternative to uniformization theorems in descriptive set theory (see e.g.~Kechris \cite{kechris1995classical}, Molchanov \cite{molchanov2005theory}, or Rockafellar and Wets   \cite{rockafellar02}), whenever one restricts attention to almost everywhere Borel selections.  
For example, we show that a maximum theorem in a conditional model is equivalent to a maximum theorem for normal integrands.  
In particular, compact subsets of a Euclidean space in a conditional model are uniquely related to compact-valued maps.  

There is a practical interest in such model-theoretic results since the transfer principle is a highly efficient tool which replaces the tedious 
work of proving by hand conditional versions of classical results which are relevant in applications.  Existing areas of application include probability (Jamneshan et al.~\cite{jamneshan2017measures}), 
mathematical economics (Backhoff and Horst \cite{backhoff2016conditional}, Bielecki et al.~\cite{bielecki2016dynamic}, Cheridito et al.~\cite{cheridito2016equilibrium}, Drapeau and Jamneshan \cite{drapeau2016conditional}, Filipovic et al.~\cite{filipovic2012approaches}, Frittelli and Maggis \cite{frittelli2011conditional,frittelli2014complete}, Hansen and Richard \cite{hansen1987role}, and Kabanov and Stricker \cite{kabanov2001teachers}), 
stochastic optimization (Cheridito and Hu \cite{cheridito2011optimal}, Cheridito and Stadje \cite{cheridito2013bsdes}, and Jamneshan et al.~\cite{jamneshan2017parameter}), set-valued analysis and measurable selection theory (Jamneshan et al.~\cite{jamneshan2017parameter} and Jamneshan and Zapata \cite{jamneshan2017compact}), Lebesgue-Bochner spaces and vector duality (Drapeau et al.~\cite{drapeau2017fenchel} and Grad and Jamneshan \cite{grad2018perturbational}), and probabilistic analysis (Guo \cite{guo2010relations,guo2013homo}, Guo et al.~\cite{guo2009separation}, Guo and Zhang \cite{guo2011neumann}, and Jamneshan and Zapata \cite{jamneshan2017compact}).

Conditional set theory (Drapeau et al.~\cite{drapeau2016algebra}), which is used to build a `conditional' model of ACA$_{0}$,   
is conceptually closely related to Boolean-valued models and topoi of sheaves, see Jamneshan \cite{jamneshan2014sheaves}.  
In fact, a conditional model of ACA$_0$ is a Boolean-valued model of ACA$_0$ by changing to the measure algebra associated to an underlying base space $(\Omega,\mathcal{F},\mu)$.  
The constructive approach of second-order arithmetic permits us to explore the semantics in a Boolean-valued model, and through this understanding build a relationship to structures and theorems in a standard model, and thus facilitate applications.    

In Aviles and Zapata \cite{aviles2018boolean}, a categorical equivalence between conditional sets and a certain class of Boolean-valued sets is discussed. 
The correspondence is based on ZFC set theory rather than second-order arithmetic. We believe that our transfer principle is of independent interest, as (1) it is an explicitly formulated theorem; (2) second-order arithmetic allows a more direct and convenient modeling of relevant mathematical notions than set theory; for example, natural numbers and real numbers are treated as primitive objects and not as complicated sets; (3) it is not argued in \cite{aviles2018boolean} how one can deduce a transfer principle from a categorical equivalence and (4) the correspondence in \cite[Theorem 3.1]{aviles2018boolean} excludes local subsets which are necessary to prove a transfer principle as shown in the proof of theorem \ref{th:model} below.  

The remainder of this paper is organized as follows. 
In section \ref{s2}, we set up the stage and collect relevant notions from conditional set theory. 
In section \ref{s3}, we prove a transfer principle for ACA$_0$. 
In section \ref{s4}, we discuss consequences and the applicability of the transfer principle.  
We conclude by a discussion of potential extensions of the transfer principle in section \ref{s5}.

\section{Preliminaries}\label{s2}

In the language $L_2$ of second-order arithmetic (see Simpson \cite{simpson2009subsystems} for an introduction), we distinguish between number variables, traditionally written in lower-case Latin letters $x,y,z,...$ and set variables, usually written as upper-case Latin letters $X,Y,Z,...$. 
Moreover, we have two constant symbols $0$ and $1$, two binary function symbols $+$ and $\cdot$ and a binary relation symbol $<$.  
Between numbers and sets exists an element relation $\in$. 
The first-order terms are number variables, constant symbols and expressions of the form $t_1+t_2$ and $t_1\cdot t_2$ for first-oder terms $t_1$ and $t_2$.  
The atomic formulas are $t_1=t_2$, $t_1<t_2$ and $t_1\in X$ where $t_1$ and $t_2$ are first-order terms and $X$ is a set variable. 
The remaining formulas are obtained from atomic formulas by the use of propositional connectives and number and set quantifiers.  

The axiomatic system of arithmetical comprehension (ACA$_{0}$) consists of the axioms for discretely ordered semirings, together with the second-order induction principle and
an axiomatic scheme saying that for any second-order formula $\phi$ containing only first-order quantifiers there is a set of all natural numbers of which $\phi$ holds, see \cite{simpson2009subsystems}. 
By the cumulative results of reverse mathematics, ACA$_{0}$ is a sufficient axiomatic basis for a great number of theorems from classical mathematics; many examples can be found in \cite{simpson2009subsystems}.

Throughout fix a $\sigma$-finite measure space $(\Omega,\mathcal{F},\mu)$, and let $\mathcal{N}$ denote its $\sigma$-ideal of null sets. 
We always identify $A,B\in\mathcal{F}$ whenever $A\Delta B\in \mathcal{N}$where $\Delta$ denotes symmetric difference.  
The resulting quotient Boolean algebra has the following relevant properties: 
\begin{itemize}
\item \emph{Completeness}: Any family in $\mathcal{F}$ has a union and an intersection in $\mathcal{F}$; 
\item \emph{Countable chain condition}: Any pairwise disjoint family in $\mathcal{F}$ is at most countable;  
\end{itemize}
see e.g.~Givant and Halmos ~\cite[Chapter 31]{halmos09} for a reference. 
We will always identify two functions $x$ and $y$ on $\Omega$ with the same codomain if $\{\omega : x(\omega)\neq y(\omega)\}\in\mathcal{N}$.  
For a function $x$ on $\Omega$ and $A\in \F$, we write $x|A$ for the restriction of $x$ to $A$. 
For a Polish space $E$, let $L^0(E)$ denote the space of Borel functions $x\colon \Omega\to E$. 
In particular, $L^0(\mathbb{R})$ denotes the space of real-valued measurable functions, and $L^0(\mathbb{N})$, $L^0(\mathbb{Z})$ and $L^0(\mathbb{Q})$ denote its subsets of functions with values in $\mathbb{N}=\{0,1,\ldots\}$, 
the integers $\mathbb{Z}$ and the rational numbers $\mathbb{Q}$ respectively.   
Recall that $L^0(\mathbb{R})$ is a Dedekind complete Riesz lattice where addition, multiplication and order are defined pointwise, see e.g.~Fremlin \cite{fremlin1974topological} for a reference. 
All inequalities between real-valued measurable functions shall be  understood in the almost everywhere sense. 
By an abuse of language, we also denote by $0$ and $1$ the functions with constant values $0$ and $1$ respectively.     
For a measurable partition $(A_k)$ and  a countable family $(x_k)$ in $L^0(E)$ for some Polish space $E$, we write $\sum_k x_k|A_k$ for the unique element $x\in L^0(E)$ such that $x|A_k=x_k|A_k$ for all $k$.

We introduce the conditional power set of $L^0(\mathbb{N})$, see Drapeau et al.~\cite{drapeau2016algebra} for an introduction to conditional set theory. 
\begin{definition} 
A set $N \subset L^0(\mathbb{N})$ is said to be \emph{stable under countable concatenations}, or \emph{stable} for short,  if 
it is not empty and 
$\sum_k n_k|A_k\in N$ for all measurable partitions $(A_k)$ and every countable family $(n_k)$ in $N$. 
The \emph{conditional power set} of $L^0(\N)$ is the collection
\[
\mathcal{P}:=\{N|A\colon N \subset L^0(\mathbb{N}) \text{ is stable, and  }A\in \F\} 
\]
where $N|A:=\{n|A\colon n\in N\}$.  
We write $N|\emptyset=\{\ast\}$.  

We define a concatenation of a countable family $(N_k|B_k)$ in $\mathcal{P}$ and a measurable partition $(A_k)$ by 
\[
\sum_k (N_k|B_k)|A_k:=\Bigg(\sum_k N_k|A_k\Bigg)|\cup_k (A_k\cap B_k)
\]
where 
\[
\sum_k N_k|A_k:=\Bigg\{\sum_k n_k|A_k\colon n_k\in N_k \text{ for each }k\Bigg\} 
\]
\end{definition}
Let us define conditional intersection and conditional complement. 
We follow the presentation in Jamneshan et al.~\cite[Section 2]{jamneshan2017measures}.  
Let $N|A,M|B\in \mathcal{P}$.  
The \emph{conditional intersection} of $N|A$ and $M|B$ is defined as 
\begin{equation}\label{eq:intersection}
N|A\sqcap M|B:=V|C 
\end{equation}
where 
\begin{align*}
C&:= \cup\{C^\prime\in \F\colon C^\prime\subset A\cap B, N|C^\prime\cap M|C^\prime\neq \emptyset\}\\
V&:=\{n\in L^0(\N)\colon n|C\in N|C\cap M|C\} 
\end{align*}
The \emph{conditional complement} of $N|A$ is the conditional subset
\begin{equation}\label{eq:complement}
(N|A)^\sqsubset:=W|D
\end{equation}
where 
\begin{align*}
D&:= \cup\{D^\prime\in \F\colon \exists\, n\in L^0(\N) \text{ such that } n|E\not\in N|E \text{ for all } E\subset D^\prime \cap A \text{ with } \mu(E)>0\}\\ 
W&:=\{n\in L^0(\N)\colon n|E\not\in N|E \text{ for all } E\subset D \cap A \text{ with } \mu(E)>0\}
\end{align*}
By applying an exhaustion argument, it can be derived from stability that $C$ and $D$ are attained, and it can 
be checked that $V$ and $W$ are stable sets as well.  We conclude that the conditional intersection and conditional complement are well defined,  see e.g.~\cite{jamneshan2017measures} for a complete argument.  

We introduce a conditional element relation between $L^0(\N)$ and its conditional power set $\mathcal{P}$ which is a new ingredient in conditional set theory.    
\begin{definition}\label{d:elementrelation}
The \emph{conditional element relation} is the function 
\begin{align*}\label{eq:element}
i &\colon L^0(\N) \times \mathcal{P}\to \mathcal{F}\\
i(n,N|A)&:=\cup\{A^\prime\in\F\colon A^\prime\subset A, n|A^\prime\in N|A^\prime\}
\end{align*}  
\end{definition}
The union of $\mathcal{G}:=\{A^\prime\in\F\colon A^\prime\subset A, n|A^\prime\in N|A^\prime\}$ is attained.    
Indeed, from Givant and Halmos, \cite[Section 30, Lemma 1]{halmos09} we know that there exists a countable family $(B_k)$ in $\mathcal{G}$ such that $A_\ast=\cup_k B_k=\cup \mathcal{G}$.  
Let $A_1=B_1$ and $A_k=B_k\cap (A_1\cup A_2\cup\ldots \cup A_{k-1})^c$ for $k\geq 2$ which defines a measurable partition of $A_\ast$.    
Since $A_k\in \mathcal{G}$ for all $k$, the claim follows from stability of $N$. 
Further, it can be directly checked that
\begin{equation}\label{eq1}
i\Bigg(\sum_k n_k|A_k, N|A\Bigg) = \cup_{k} (i(n_k,N|A) \cap A_k) 
\end{equation}

\section{A Transfer Principle for ACA$_0$}\label{s3}

Our aim is to prove that every consequence of ACA$_{0}$ holds in the structure $\mathcal{S}:=(L^0(\N),\mathcal{P},+,\cdot,0,1,<, i)$ with truth value $\Omega$, i.e.~that $\mathcal{S}$ is a `conditional model' of ACA$_0$. 
We start by explaining the evaluation of terms and formulas in $\mathcal{S}$.

Let $V_1$ denote the collection of all number variables and let $V_2$ denote the collection of all set variables.  
Let $\beta$ be a function with domain $V_1\cup V_2$ such that $\beta:V^1\rightarrow L^0(\mathbb{N})$ and $\beta:V^{2}\rightarrow\mathcal{P}$, and let $t$ be a first-order term. 
Such a $\beta$ is called a `conditional assignment'.
Then $[t]^{\beta}$, the $\beta$-evaluation of $t$, is defined recursively as 
\begin{itemize}
\item $[a]^{\beta}=a$ for $a\in\{0,1\}$, 
\item $[x]^{\beta}=\beta(x)$ for $x\in V_1$, 
\item $[X]^{\beta}=\beta(X)$ for $X\in V_2$, 
\item $[t_{0}\circ t_{1}]^{ \beta}=[t_{0}]^{ \beta}\circ[t_{1}]^{ \beta}$ for $\circ\in\{+,\cdot\}$. 
\end{itemize}
For first-order terms $t_0$ and $t_1$ and a second-order variable $X$, the \emph{conditional $\beta$-evaluation} of atomic formulas is defined by
\begin{itemize}
\item $[t_0=t_1]^{ \beta}=\{\omega\colon  [t_0]^{\beta}(\omega)=[t_1]^{\beta}(\omega)\}$,  
\item $[t_0<t_1]^{ \beta}=\{\omega\colon  [t_0]^{\beta}(\omega)<[t_1]^{\beta}(\omega)\}$, 
\item $[t_1\in X]^{ \beta}=i([t_1]^{\beta}, \beta(X))$.    
\end{itemize}
The conditional $\beta$-evaluation of composite and quantified formulas is defined by 
\begin{itemize}
\item $[\phi\wedge \psi]^{ \beta}=[\phi]^{ \beta}\cap [\psi]^{ \beta}$, 
\item $[\neg \phi ]^{ \beta}=([\phi]^{ \beta})^c$, 
\item $[\exists x \phi(x)]^{ \beta}=\cup_{n\in L^0(\mathbb{N})} [\phi(x)]^{ \beta[\frac{n}{x}]}$,  
\item $[\exists X \phi(X)]^{ \beta}=\cup_{N|A \in \mathcal{P}} [\phi(X)]^{ \beta[\frac{N|A}{X}]}$. 
\end{itemize}
The remaining composite and quantification formulas are defined in the obvious way. 

We have the following maximum principle, also known from Boolean-valued models, see e.g.~\cite[Chapter 1]{bell2005set}. 
\begin{proposition}\label{p1}
Let $\beta$ be a conditional assignment and let $\phi$ and $\psi$ be formulas in $L_2$. 
Then there exist $n\in L^0(\mathbb{N})$ and $N|A\in \mathcal{P}$ such that 
\begin{align*}
[\exists x \phi(x)]^{ \beta}&=[\phi(x)]^{ \beta[\frac{n}{x}]} \\
[\exists X \psi(X)]^{ \beta}&=[\psi(X)]^{ \beta[\frac{N|A}{X}]}
\end{align*}
\end{proposition}
\begin{proof}
We may assume that $[\exists x \phi(x)]^{ \beta}=\Omega$. 
We find a countable family $(A_k=[\phi(x)]^{ \beta[\frac{n_k}{x}]})$ such that $\cup_{k} A_k=[\exists x \phi(x)]^{ \beta}$. 
Form a measurable partition from $(A_k)$, still denoted by $(A_k)$. 
Put $n=\sum_k n_k|A_k$. 
Then it holds that 
\[
[\phi(x)]^{ \beta[\frac{n}{x}]}\supset [\phi(x)]^{ \beta[\frac{n_k}{x}]}\cap [n_k=n]^\beta
\]
for all $k$ which implies $[\phi(x)]^{ \beta[\frac{n}{x}]}=[\exists x \phi(x)]^{ \beta}$.  
The second claim can be shown analogously by using concatenations in $\mathcal{P}$. 
\end{proof}

We will now adapt the usual notion of the correctness of a sequent to the conditional context.

\begin{definition}
If $\Gamma$ and $\Delta$ are sets of second-order formulas, then $\Gamma\rightarrow \Delta$ is called a \emph{sequent}.  
The \emph{conditional validity} of a sequent $\Gamma\rightarrow \Delta$ with respect to a conditional assignment $\beta$ is defined by 
\[
 (\cup_{\psi \in \Gamma} [\neg \psi]^{ \beta}) \cup (\cup_{\psi \in \Delta} [\psi]^{ \beta}) 
\]
and $\Gamma\rightarrow\Delta$ is said to be \emph{correct} if and only if
\[
(\cup_{\psi \in \Gamma} [\neg \psi]^{ \beta}) \cup (\cup_{\psi \in \Delta} [\psi]^{ \beta})=\Omega 
\]
for all assignments $\beta$. 
An inference rule $R$ is a pair consisting of a finite sequence $(\Gamma_{1},\Gamma_{2},\ldots,\Gamma_{n})$ of sequents and a single 
sequent $\Gamma$ written as $\frac{\Gamma_{1},\Gamma_{2},\ldots,\Gamma_{n}}{\Gamma}$, and 
it is said to be correct if and only if the correctness of $\Gamma$ follows from the correctness of $\Gamma_1,\Gamma_{2},\ldots,\Gamma_n$. 
\end{definition}
We want to apply the inference rules for second-order logic given by Takeuti in \cite[p.~9-10 and p.~135-136]{takeuti2013proof}. 
By Boolean arithmetic, one can directly check that all structural and logical rules are correct.     
We illustrate this for the first weakening rule 
\begin{equation}\label{eq4}
\frac{\Gamma\rightarrow\Delta}{\phi,\Gamma\rightarrow\Delta}
\end{equation}
For each assignment $\beta$, one has 
\begin{align*}
(\cup_{\psi \in \Gamma} [\neg \phi]^{ \beta}) \cup (\cup_{\psi \in \Delta} [\psi]^{ \beta}) \subset  (\cup_{\psi \in \Gamma\cup \{\phi\}} [\neg \phi]^{ \beta}) \cup (\cup_{\psi \in \Delta} [\psi]^{ \beta})
\end{align*}
which leads to the correctness of \eqref{eq4}.  
By Proposition \ref{p1}, the left universal quantification rule 
\[
\frac{\phi(t),\Gamma\rightarrow\Delta}{\forall{x}\phi(x),\Gamma\rightarrow\Delta}
\]
where $t$ is a term, is correct.  
As for the correctness of the right universal quantification rule 
\[
\frac{\Gamma\rightarrow\Delta,\phi(y)}{\Gamma\rightarrow\Delta,\forall{x}\phi(x)}
\]
where $y$ does not occur freely in $\Gamma\rightarrow\Delta,\forall{x}\phi(x)$, assume 
\begin{equation}\label{eq2}
\cup_{\psi\in \Gamma} [\neg \psi]^{\beta}\cup \cup_{\psi\in \Delta} [\psi]^{\beta} \cup [\phi(y)]^{\beta}=\Omega 
\end{equation}
for all conditional assignments and let $\beta$ be an arbitrary conditional assignment. 
We want to show that 
\begin{equation}\label{eq5}
\cup_{\psi\in \Gamma} [\neg \phi]^{\beta} \cup \cup_{\psi\in \Delta} [\psi]^{\beta}\cup \cap_{n\in L^0(\mathbb{N})}[\phi]^{\beta[\frac{n}{x}]}=\Omega
\end{equation}
By Proposition \ref{p1}, there exists $n_\ast\in L^0(\mathbb{N})$ such that $\cap_{n\in L^0(\mathbb{N})}[\phi(x)]^{\beta[\frac{n}{x}]}=[\phi(x)]^{\beta[\frac{n_\ast}{x}]}$. 
Choose a $y$-variant $\beta^{\prime}$ of $\beta$ that maps $y$ to $n_\ast$. 
Then \eqref{eq5} follows from \eqref{eq2} since $y$ does not occur freely in $\Gamma\rightarrow\Delta,\forall{x}\phi(x)$, so that the first two sets of the union \eqref{eq5} remain unchanged.  
Analogously, one shows the first-order existential quantification inference rules. 
The second-order quantifier inference rules are only relevant for second-order variables, as predicate constants do not appear in our language.
The inference rules for second-order quantification can hence be proved analogously to the corresponding first-order rules.    
Thus, we obtain:

\begin{lemma}{\label{deduction}}
If $\frac{\Gamma_{1},...,\Gamma_{n}}{\Gamma}$ is any deduction rule of second-order sequent calculus and $\Gamma_{1},...,\Gamma_{n}$ are correct, then $\Gamma$ is correct.
In particular, if all elements of $\Gamma$ hold in $\mathcal{S}$ with truth value $\Omega$ and $\Gamma\rightarrow\phi$ is derivable with the rules of sequent calculus,
then $\phi$ holds in $\mathcal{S}$ with truth value $\Omega$.
\end{lemma}

\begin{theorem}\label{th:model}
All the axioms of ACA$_0$ attain the value $\Omega$ in the structure $\mathcal{S}$ for all conditional assignments. 
\end{theorem}
\begin{proof}
Let $\beta$ be an arbitrary conditional assignment. 
The verification of the basic axioms (see Simpson \cite[p.~4]{simpson2009subsystems}) is immediate from the definitions.  
For the sake of completeness, we provide the elementary arguments below. 
\begin{itemize}
\item $[\neg (x+1=0)]^{ \beta}=([x+1=0]^{ \beta})^c=\{\omega\colon \beta(x)(\omega)+1=0\}^c=\emptyset^c=\Omega$. 
\item Since $\{\omega\colon \beta(x)(\omega)+1=\beta(y)(\omega)+1\}=\{\omega\colon \beta(x)(\omega)=\beta(y)(\omega)\}$, it follows from $[x+1=y+1]^{ \beta}=\Omega$ that $[x=y]^{ \beta}=\Omega$.
\item Similarly, one can verify that $[x+0=x]^{ \beta}=\Omega$, $[x+(y+1)=(x+y)+1]^{ \beta}=\Omega$, $[x\cdot 0=0]^{ \beta}=\Omega$, $[x\cdot (y+1)=(x\cdot y)+x]^{ \beta}=\Omega$, and $[\neg(x<0)]^{\beta}=\Omega$. 
\item $[x<y+1]^{\beta}=\Omega$ is equivalent to $\{\omega\colon \beta(x)(\omega)=\beta(y)(\omega)\}\cup \{\omega\colon \beta(x)(\omega)<\beta(y)(\omega)\}=\Omega$, which means that $[(x<y)\vee (x=y)]^{ \beta}=\Omega$. 
\end{itemize}
As for the second-order induction scheme, we have to verify that 
\[
[(0\in X \wedge \forall x(x\in X \rightarrow  x+1\in X))\rightarrow \forall x(x\in X)]^{ \beta}=\Omega
\]
By conditionally evaluating the previous formula and rewriting it by using Boolean arithmetic, we must verify that $A\subset B$, where 
\begin{align*}
A&:=i(0,\beta(X))\cap \cap_{n\in L^0(\mathbb{N})} (([i(x,X)]^{\beta[\frac{n}{x}]})^c \cup [i(x+1,X)]^{\beta[\frac{n}{x}]})\\
B&:=\cap_{n\in L^0(\mathbb{N})} [i(x,X)]^{\beta[\frac{n}{x}]}
\end{align*}
But this is immediate from the stability of $\beta(X)$.  

Finally, we verify the arithmetical comprehension scheme, that is, we want to show that 
\[
[\exists X \forall x (x\in X \leftrightarrow \phi(x))]^\beta=\Omega 
\]
for any arithmetical\footnote{Recall that a formula of $L_2$ is said to be \emph{arithmetical} if it contains no set quantifiers, see e.g. Simpson~\cite{simpson2009subsystems}.} formula $\phi(x)$ in which $X$ does not occur freely. 
By Proposition \ref{p1}, $A_\phi:=\cup_{n\in L^0(\N)} [\phi(x)]^{\beta[\frac{n}{x}]}$ is attained.  
Suppose for a moment that 
\[
N_\phi:=\{n\in L^0(\N)\colon [\phi(x)]^{\beta[\frac{n}{x}]}=A_\phi\}
\]
satisfies stability.  
Then 
\begin{equation}\label{eq3}
[\phi(x)]^{\beta[\frac{n}{x}]}=i(n,N_\phi|A_\phi)
\end{equation}
for all $n\in L^0(\N)$.  
Indeed, for $n_0\in L^0(\N)$, let $n_1\in N_\phi$ be such that $n_0|B=n_1|B$ where $B=i(n_0,N_\phi|A_\phi)$, and put $n_2=n_0|B+n_1|B^c$.  
By stability of $N_\phi$, we have 
\[
[\phi(x)]^{\beta[\frac{n_0}{x}]}=[\phi(x)]^{\beta[\frac{n_2}{x}]}\cap [n_2=n_0]^\beta=A_\phi\cap B=B
\] 
By Boolean arithmetic, it follows from \eqref{eq3} that 
\[
[\forall x (x\in X \leftrightarrow \phi(x))]^{\beta[\frac{Y_\phi|A_\phi}{X}]}=\Omega  
\]
which proves comprehension.  
Thus it remains to verify that $N_\phi$ is stable under concatenations for all formulas $\phi$, which we will prove by an induction on arithmetical formulas.  
First, since addition and multiplication commute with 
concatenations\footnote{That is, $\sum_k n_k|A_k\circ \sum_k m_k|B_k=\sum_{k,h} (n_k\circ m_h)|A_k\cap B_h$ for $\circ\in\{+,\cdot\}$.}, for any first-order term $t=t(x)$, by an induction on terms, one has
\[
[t(x)]^{\beta[\frac{\sum_k n_k|A_k}{x}]}=\sum_k [t(x)]^{\beta[\frac{n_k}{x}]}|A_k
\]
for all measurable partitions $(A_k)$ and every countable family $(n_k)$ in $L^0(\N)$.  
Since also order and concatenations commute and due to \eqref{eq1}, $N_\phi$ is stable for all atomic formulas $\phi$. 

Let $\phi$ and $\psi$ be two arithmetical formulas such that $N_\phi|A_\phi, N_\psi|A_\psi\in \mathcal{P}$. 
Then we have 
\begin{align*}
[(\phi\wedge \psi)(x)]^{\beta[\frac{n}{x}]}&=[\phi(x)]^{\beta[\frac{n}{x}]} \cap [\psi(x)]^{\beta[\frac{n}{x}]}\\
&=i(n,N_\phi|A_\phi)\cap i(n,N_\psi|A_\psi) \\
&=i(n,N_\phi|A_\phi \sqcap N_\psi|A_\psi)
\end{align*}
Moreover, for a negation one obtains 
\[
[\neg \phi(x)]^{\beta[\frac{n}{x}]}=i(n,N_\phi|A_\phi)^c=i(n,(N_\phi|A_\phi)^\sqsubset)
\]
Finally, let $\theta(x,y)$ be arithmetical and $\phi(y)=\exists x \theta(x,y)$.  
Clearly, $A_\phi=A_\theta$.  
By the established, we can already define a pairing function, product of stable sets and their projections, see e.g. Simpson~\cite[p.~66-69]{simpson2009subsystems}.  
The pairing function $(i,j)\mapsto (i+j)^2 +i$ underlying the definition of a product commutes with concatenations since this is the case for addition and multiplication.  
Whence if 
\[
N_\theta=\{(n,m)\in L^0(\N)^2\colon [\theta(x,y)]^{\beta[\frac{n}{x},\frac{m}{y}]}=A_\theta\}
\]
is stable, then its projection to the first coordinate is stable, and by definition equal to $N_\phi$. 
\end{proof}

We now state the main theorem of this section. 

\begin{theorem}\label{transfer}
If $\phi$ is a consequence of ACA$_{0}$, then $\phi$ holds in $\mathcal{S}$ with truth value $\Omega$.
\end{theorem}
\begin{proof}
 By theorem \ref{th:model}, all axioms of ACA$_{0}$ have truth value $\Omega$ in $\mathcal{S}$.
 By Lemma \ref{deduction}, truth value $\Omega$ is preserved under the rules of second-order sequent calculus. 
 Thus, the theorem follows.
\end{proof}

\section{Harvesting the fruits}\label{s4} 

A first step towards applications of the model-theoretic results in this article is to investigate connections between a conditional and a standard model. 
An aim of this section is to understand consequences of the transfer principle and interpret them properly in a standard setting. 
For example, the structure $L^0(\N)$ are the natural numbers in the 'conditional model', while interpreted in a standard model, it is the space of $\N$-valued measurable functions on a $\sigma$-finite measure space.  
We shall discover connections to different parts of analysis such as set-valued analysis and measurable selection theory, vector-valued analysis and probabilistic analysis with applications to areas such as stochastic control, vector optimization and mathematical economics. 
An important observation is that application of theorems in a 'conditional model' preserves measurability systematically by construction.   
Moreover, we will demonstrate that a transfer principle unifies a conditional version of classical theorems which were proved by hand previously.  Finally, 
new examples of theorems and applications shall underline the usefulness of a transfer principle.

\subsection{Set theory}

We begin by collecting basic set-theoretical vocabulary in a conditional setting.
The  \emph{inclusion relation} (see Simpson \cite[Definition II.3.1]{simpson2009subsystems}) on the conditional power set $\mathcal{P}$ can be interpreted by the relation 
\begin{equation}\label{eq:inclusion}
N|A\sqsubseteq M|B \quad \text{ if and only if }\quad A\subseteq B\;\text{ and } \;N|A\subseteq M|A. 
\end{equation}
This relation coincides with the conditional inclusion relation in conditional set theory, see Drapeau et al.~\cite[Definition 2.8]{drapeau2016algebra} for the abstract formulation in the context of an arbitrary complete Boolean algebra, and see Jamneshan et al.~\cite[Definition 2.5]{jamneshan2017measures} for a formulation in the context of an associated measure algebra. 

The \emph{product} of two sets $N|A,M|B\in\mathcal{P}$ is the set 
\begin{equation}\label{eq:product}
N\times M|A\cap B
\end{equation}
by interpreting \cite[Definition II.3.1]{simpson2009subsystems} in $\mathcal{S}$.   
This definition extends the definition of a conditional product in conditional set theory \cite[Definition  2.14]{drapeau2016algebra}.  
Indeed, a conditional product was introduced in \cite{drapeau2016algebra} only for sets of the form $N=N|\Omega$ (if we consider the setting of the present paper).  
It was realized in \cite{jamneshan2017measures} that the definition given in \cite{drapeau2016algebra} does not suffice for proving a conditional version of Fubini's theorem.  
The extended definition then used in \cite[Section 5.1]{jamneshan2017measures} coincides with \eqref{eq:product}. 

Following \cite[Definition II.3.1]{simpson2009subsystems}, a \emph{function} $f\colon N|A\to M|B$ can be interpreted as a subset $W|C \sqsubseteq N\times M|A\cap B$ such that for each $n|C\in N|C$ there is a unique $m|C\in M|C$ such that $(n,m)|C\in W|C$.  
This definition extends the one of a conditional function in \cite[Definition 2.17]{drapeau2016algebra} (similarly as the definition of a product \eqref{eq:product} extends the definition of a conditional product in \cite{drapeau2016algebra}). 
The concept of a conditional function is a significant abstraction of earlier concepts named a regular function in Detlefsen and Scandolo \cite[Definition 4, Proposition 1]{detlefsen2005conditional}, a stable function in Cheridito et al.~\cite[Definition 4.2]{cheridito2015conditional}, and a local function in Filipovic et al.~\cite[Definition 3.1]{filipovic2009separation}. 
The characteristic property of a conditional function is that evaluation of its values commutes with concatenations.

From \cite[Lemma II.2.1]{simpson2009subsystems} we know that $(L^{0}(\mathbb{N}),+,\cdot,0,1,<)$ is a commutative ordered semiring with cancellation.  
For instance, the \emph{totality property} $m<n \vee m=n \vee n<m$ reads as the statement that for each pair $n,m\in L^0(\N)$ there is a partition $(\{m<n\},\{m=n\}, \{n<m\})$. 
This interpretation coincides with the definition of conditionally total in \cite[Definition 2.15]{drapeau2016algebra}. 

A set $N|A$ is said to be \emph{finite} whenever there exists $k\in L^0(\N)$ such that $n|A<k|A$ for all $n|A\in N|A$, see \cite[p.~67]{simpson2009subsystems}, which coincides with the notion of conditionally finite, see \cite[Definition 2.23]{drapeau2016algebra}.  In general, a finite set in the model  $ \mathcal{S} $ is not finite in a standard sense, it might even be uncountable from the latter perspective. However, for every finite set $N|A$ there exist a partition $(A_k)$ and a countable family $(N_k)$ of finite subsets of $\mathbb{N}$ such that $N|A=(\sum_k N_k|A_k)|A$, where we identified $N_k$ with the set of all measurable $n\colon \Omega\to N_k$, cf.~\cite[Lemma 2.22]{drapeau2016algebra}. 
The \emph{collection of all finite sequences} of length $n\in L^0(\N)$ is the set of functions $\{1\leq m\leq n\}\to L^0(\N)$, see \cite[Definition II.3.3]{simpson2009subsystems}, which is precisely the construction given in the paragraph after \cite[Definition 2.20]{drapeau2016algebra}.  

A \emph{sequence} in a set $N$ is a function $f\colon L^0(\N)\to N$, compare with the definition of a conditional sequence in \cite[Example 2.2.1]{drapeau2016algebra}.  
From a standard point of view such a sequence is a net which is parametrized by $L^0(\N)$ and commutes with concatenations,  cf. Drapeau et al.~\cite[Section 2]{drapeau2017fenchel}. 

Let $(N_k|A_k)$ be a sequence in the power set $\mathcal{P}$.  
Using the element relation (see Definition \ref{d:elementrelation}), we see that the \emph{intersection} of $(N_k|A_k)$ can be identified with  
\begin{equation}\label{eq:inter}
\sqcap_k N_k|A_k:=V|B
\end{equation}
where 
\begin{align*}
B&:= \cup\{B^\prime\in \F\colon B^\prime\subset \cap_k A_k, \; \cap_k (N_k|B^\prime)\neq \emptyset\}\\
V&:=\{n\in L^0(\N)\colon n|B\in \cap_k (N_k|B)\}
\end{align*}
and that the  \emph{union} of $(N_k|A_k)$ can be identified with 
\begin{equation}\label{eq:union}
\sqcup_k N_k|A_k:=W|C
\end{equation}
where 
\begin{align*}
C&:= \cup_k A_k \\
W&:=\{n\in L^0(\N)\colon \cup_k (i(n,N_k)\cap A_k)= C\}
\end{align*}
The \emph{complement} of a set was defined in \eqref{eq:complement}.  
We have seen that the set operations in $ \mathcal{S} $ recover the conditional set operations \cite[p.~567]{drapeau2016algebra}.  
In \cite[Corollary 2.10]{drapeau2016algebra}, it was proved that the conditional power set has the structure of a complete Boolean algebra, which is a fundamental result for basic constructions in conditional topology \cite[Section 3]{drapeau2016algebra} and conditional measure theory \cite{jamneshan2017measures}. 
The transfer principle for ACA$_0$ implies a weaker statement, namely that the power set $\mathcal{P}$ has the structure of a $\sigma$-complete Boolean algebra, see \cite{simpson2009subsystems}.

\subsection{Real analysis and linear algebra}\label{sec:analina} 

A detailed construction of the conditional real numbers and their conditional algebraic, order and topological properties are developed in Jamneshan \cite[Chapter 5]{jamneshan2014theory} where all properties are proved from conditional set theory. 
We will argue that most of these properties and related results are consequences of the transfer principle.  

In the 'conditional model' $ \mathcal{S} $, the integers (see Simpson \cite[p.~73]{simpson2009subsystems}) are the space $L^0(\mathbb{Z})$ of measurable integer-valued functions.
From \cite[Theorem II.4.1]{simpson2009subsystems} we know that $L^0(\mathbb{Z})$ is an Euclidean ordered integral domain.  
Similarly, the space $L^0(\mathbb{Q})$ of measurable rational-valued functions are the rational numbers in $\mathcal{S}$. 
By \cite[Theorem II.4.2]{simpson2009subsystems}, $L^0(\mathbb{Q})$ is an ordered field; compare this with the  conditionally ordered field of conditional rational numbers in \cite[Example 4.2.1]{drapeau2016algebra}.  
The topological completion of $L^0(\mathbb{Q})$ inside $\mathcal{S}$ are the real numbers which by a standard approximation argument can be identified with the space $L^0(\mathbb{R})$ of real-valued measurable functions, see \cite[Definition II.4.4]{simpson2009subsystems}.  
Now \cite[Theorem II.4.5]{simpson2009subsystems} implies that $L^0(\mathbb{R})$ is an Archimedean ordered field, cf.~\cite[Lemma 5.2.12 and Theorem 5.2.7]{jamneshan2014theory}.

Notice that the absolute value of real numbers maps to $L^0_+(\mathbb{R}):=\{x\in L^0(\mathbb{R})\colon x\geq 0\}$. 
Let $L^0_{++}(\mathbb{Q}):=\{x\in L^0(\mathbb{Q})\colon x> 0\}$. 
 An \emph{open ball} in $L^0(\mathbb{R})$ is a set 
 \[
 B(q,r):=\{x\in L^0(\mathbb{R})\colon |x-q| <r\}
 \]
where $q\in L^0(\mathbb{Q})$ and $r\in L^0_{++}(\mathbb{Q})$, see \cite[p.~81]{simpson2009subsystems}.    
An \emph{open set} $O$  is the union  of a sequence of open balls $\sqcup_k B(q_k,r_k)$,  see \cite[Definition II.5.6]{simpson2009subsystems}.  
A \emph{closed set} is the complement of an open set, see \cite[Definition II.5.12]{simpson2009subsystems}.  
For example, for any $r\in L^0_{++}(\mathbb{Q})$, the set $\{x\in L^0(\mathbb{R})\colon |x| \leq r \}$ is closed, but the standard complement of $B(0,r)$ is in general not.  
However, it can be easily verified that every  $\mathcal{S}$-closed set is sequentially closed, i.e.~it  contains the limit of any almost everywhere converging standard sequences.  

\begin{remark}\label{r:realline}
The Euclidean topology of the real numbers $L^0(\mathbb{R})$ in $\mathcal{S}$ is the order topology in a standard model which renders pointwise addition and multiplication continuous. 
Thus, $L^0(\mathbb{R})$ becomes a topological algebra, in particular a topological module of rank $1$ over itself. 
This fact is exploited in Filipovic et al.~\cite{filipovic2009separation}, Cheridito et al.~\cite{cheridito2015conditional}, Drapeau et al.~\cite{drapeau2016algebra}, 
Jamneshan and Zapata \cite{jamneshan2017compact} to build a functional analytic discourse in $L^0(\mathbb{R})$-modules. 
We shall observe later that a functional analytic discourse in $L^0(\mathbb{R})$-modules is the reflection of classical functional analysis, albeit in a 'conditional model'.  
So far, we know that $L^0(\mathbb{R})$ is $1$-dimensional Euclidean space in the 'conditional model' $\mathcal{S}$. 
\end{remark}
Let $L^0(\mathbb{R})^n$ denote the set of finite sequences $\{1\leq m \leq n\}\to L^0(\mathbb{R})$.  
Write $n$ as $\sum_k n_k|A_k$ for $(n_k)$ in $\mathbb{N}$ and $(A_k)$ a measurable partition. 
Then one can view $L^0(\mathbb{R})^n$ as $\sum_k L^0(\mathbb{R})^{n_k}|A_k$ where each $L^0(\mathbb{R})^{n_k}$ is a standard product. 
The absolute value extends from $L^0(\mathbb{R})$ to an \emph{Euclidean norm} on $L^0(\mathbb{R})^n$ which for $x=\sum_k (x_1,\ldots,x_{n_k})|A_k$ is defined by  
\begin{equation}\label{eq:norm}
\|x\|:=\sum_k \sqrt{x_1^2+\ldots+x_{n_k}^2}\bigg|A_k
\end{equation}
Convergence in the conditional Euclidean space $L^0(\mathbb{R})^n$ can be characterized by almost everywhere convergence. 
Suppose that $n\in \mathbb{N}$ (the general case $n=\sum_k n_k|A_k$ follows by localizing the subsequent construction to each $A_k$ and gluing).   
Let $(x_k)$ be a standard sequence in $L^0(\mathbb{R})^n$ converging almost everywhere to $x$. 
Given $r\in L^0_{++}(\mathbb{Q})$, let
\[
k_r(\omega)=\inf\{k^\prime\in \N\colon \|x_m -x\|(\omega)< r(\omega) \text{ for all } m \geq k^\prime\}
\]
Notice that $k_r$ is measurable. 
Obtain from $(x_k)$ the conditional sequence $x_k:=\sum_n x_{k_n}|A_n$, $k=\sum_n k_n|A_n$. 
By construction, $\|x_k-x\|<r$ for all $k\geq k_r$.  
Conversely, it is easy to see that if $(x_k)_{k\in L^0(\N)}$ is a conditional sequence conditionally converging to $x$, 
then the standard subnet $(x_k)_{k\in \N}$ resulting by embedding $\N$ into $L^0(\N)$ via $n\mapsto n 1_\Omega$ converges to $x$ almost everywhere.

\begin{remark}\label{r:euclidean}
The $\mathcal{S}$-Euclidean norm \eqref{eq:norm} is an example of $L^0(\R)$-valued vector norms which appear in different parts of analysis such as the analysis of Lebesgue-Bochner spaces and vector integration (see e.g.~Diestel and Uhl~\cite{joseph1977vector}, Haydon et al.~\cite{haydon1991randomly}, Hyt\"onen et al.~\cite{hytoenen2016analysis}),  financial modeling (see e.g.~Hansen and Richard \cite{hansen1987role}, Filipovic et al.~\cite{filipovic2012approaches}, and Cheridito et al.~\cite{cheridito2016equilibrium}), or probabilistic analysis (see e.g.~Guo \cite{guo2010relations} and Haydon et al.~\cite{haydon1991randomly}).  
In section \ref{sec:banach}, we will continue to discuss the connection of $L^0(\R)$-valued vector normed spaces and Banach spaces in the 'conditional model' $\mathcal{S}$.   
\end{remark}
By the transfer principle, we obtain from \cite[Lemma III.2.1]{simpson2009subsystems} a Bolzano-Weierstra\ss~theorem in $\mathcal{S}$. 
Recall that a sequence $(x_k)$ in $L^0(\mathbb{R})^n$, $n\in L^0(\N)$, is said to be \emph{bounded} if there exists $r\in L^0_{++}(\mathbb{Q})$ such that $\|x_k\|<r$ for all $k\in L^0(\N)$.   
\begin{theorem}\label{t:BW}
Let $(x_k)$ be a bounded sequence in $L^0(\mathbb{R})^n$. Then $x=\limsup_k x_k$ exists.   
Moreover, there exists a subsequence $(x_{k_p})$ converging to $x$.  
\end{theorem}
The previous statement is reminiscent of  a `measurable' or `conditional' version of the Bolzano-Weierstra\ss~theorem as proved in Kabanov and Stricker \cite[Lemma 2]{kabanov2001teachers},  in F\"ollmer and Schied \cite[Lemma 1.64]{foellmer2011stochastic} and in Cheridito et al.~\cite[Theorem 3.8]{cheridito2015conditional} respectively. 
Indeed, any standard sequence $(x_k)$ in $L^0(\mathbb{R})$ can be extended to a conditional sequence (see above). 
If $(x_k)$ is bounded ($\sup_k |x_k|<\infty$), then there is a conditional subsequence $(x_{k_p})$ which converges to $x=\limsup_k x_k$ by theorem \ref{t:BW}.  
Choosing the standard subsequence $(x_{k_p})_{p\in \N}$ in the net $(x_{k_p})_{p\in L^0(\N)}$, we conclude $x_{k_p}\to x$ almost everywhere.

Compactness within ACA$_0$ is introduced in \cite[Definition III.2.3]{simpson2009subsystems} as a form of sequential compactness.  
We have the following characterization.  
\begin{proposition}\label{p:compact}
Let $W$ be a stable subset of $L^0(\R)^n$, $n=\sum_k n_k|A_k\in L^0(\mathbb{N})$. Then the following are equivalent. 
\begin{itemize}
\item[(i)] $W$ is conditionally compact. 
\item[(ii)] \emph{Heine-Borel property}: For every conditional sequence $(O_i)$ of open sets such that $W\sqsubseteq \sqcup_i O_i$ there exists a conditionally finite subsequence $(O_{i_j})$ such that $W \sqsubseteq \sqcup_j O_{i_j}$.  
\item[(iii)] $W$ can be represented as $\sum_k W_k|A_k$ where each $W_k$ is the set of almost everywhere selections of an Effros measurable\footnote{A set-valued function $X\colon\Omega\to 2^{\R^n}$ is said to be Effros measurable, if $\{\omega\colon X(\omega)\cap O\neq \emptyset\}\in \mathcal{F}$ for all open sets $O$ in $\mathbb{R}^n$.} compact-valued map in $\mathbb{R}^{n_k}$. 
\end{itemize}
\end{proposition}
\begin{proof}
The equivalence of $(i)$ and $(ii)$ follows from \cite[Theorem IV.1.5]{simpson2009subsystems} and the transfer principle \ref{transfer}. 
The equivalence of $(ii)$ and $(iii)$ is proved in Jamneshan and Zapata \cite[Section 5]{jamneshan2017compact}.  
\end{proof}

\begin{remark}
The notion of conditional compactness was introduced in \cite[Definition 3.24]{drapeau2016algebra}. 
A conditional version of the Heine-Borel theorem for general conditional metric spaces was established in \cite[Theorem 4.6]{drapeau2016algebra}. 
A simple example of a conditionally compact set in $L^0(\mathbb{R})$ are finite unions of intervals  
\[
[a,b]:=\{x\in L^0(\mathbb{R})\colon a\leq x\leq b\}
\]	
where $a<b$ in $L^0(\mathbb{R})$. 
\end{remark}
The following minimum theorem was proved in \cite[Theorem 4.4]{cheridito2015conditional} in finite dimensions, and in full generality  in Jamneshan and Zapata \cite[Theorem 5.13]{jamneshan2017compact}, and applied in e.g. Cheridito et al.~\cite{cheridito2016equilibrium} and Jamneshan et al.~\cite{jamneshan2017parameter} successfully to stochastic control.   
\begin{theorem}\label{t:max}
Let $W$ be a conditionally compact subset of $L^0(\mathbb{R})^n$, $n\in L^0(\N)$, and let $f\colon W\to L^0(\mathbb{R})$ be a conditionally lower semi-continuous function.  
Then $f$ has a minimum.  
\end{theorem}
\begin{proof}
The proof can be done in WKL$_0$ which is a subsystem of ACA$_{0}$ using the Heine-Borel covering property \cite{simpson2009subsystems}. The claim then follows from theorem \ref{transfer}.
\end{proof}

We can derive the following variant of the previous theorem in set-valued analysis an important aspect of which is the avoidance of measurable selection arguments.   
\begin{corollary}
Let $f\colon \Omega\times \mathbb{R}^n\to \mathbb{R}$ be a normal integrand\footnote{A function $f\colon \Omega\times \mathbb{R}^n\to \mathbb{R}$ is said to be a normal integrand if its graph is Effros measurable and closed-valued.}, and let $X\colon \Omega\to 2^{\mathbb{R}^n}$ be an Effros measurable compact-valued map. 
Then there exists a measurable function $x\colon \Omega\to \mathbb{R}^n$ such that $f(\omega,x(\omega))=\min_{x\in X(\omega)} f(\omega, x)$ almost everywhere. 
\end{corollary}
\begin{proof}
By proposition \ref{p:compact}, one can identify $X$ with a conditionally compact set in $L^0(\mathbb{R})^n$. 
From a result in Jamneshan et al. \cite[Section 5]{jamneshan2017parameter}, one can identify $f$ with a sequentially lower semi-continuous\footnote{Lower semi-continuity in $\mathcal{S}$ can be interpreted as sequential lower semi-continuity in a standard model, i.e.~$f(x)\leq \liminf f(x_k)$ almost everywhere whenever $(x_k)$ converges almost everywhere to $x$, see the discussion above relating convergence in the $\mathcal{S}$-Euclidean topology with almost everywhere convergence.} function $L^0(\mathbb{R})^n \to L^0(\mathbb{R})$.  
Theorem \ref{t:max} proves the claim. 
\end{proof}
A notion of an $L^0(\R)$-derivative for functions $f\colon L^0(\mathbb{R})^n\to L^0(\mathbb{R})$ is introduced in Cheridito et al. \cite[Section 7]{cheridito2015conditional}.  An interpretation of the definition of a derivative in second-order arithmetic yields the same concept which is defined below for $n=1$. 
\begin{definition}
Let $f\colon L^0(\mathbb{R})\to L^0(\mathbb{R})$ be a stable function such that $f(x_k)\to f(x)$ almost everywhere whenever $x_k\to x$ almost everywhere.  
The value 
\[
f^\prime(x)=\lim_{h\to 0}\frac{f(x-h)-f(x)}{h}
\]
is called the \emph{conditional derivative} of $f$ at $x$ if it exists.  
\end{definition}
The following conditional version of the \emph{Peano existence theorem}, derived from \cite[Theorem IV.8.1]{simpson2009subsystems}, has not been proved previously. 
It can be applied to solve random ordinary differential equations.  
\begin{theorem}
Let $a,b\in L^0(\mathbb{R})$ with $a,b>0$. Let $f\colon [-a,a]\times [-b,b]\to L^0(\mathbb{R})$ be stable sequentially continuous. 
Then the random ordinary differential equation 
\begin{align*}
\frac{dy}{dx}=f(x,y), \quad y(0)=0, 
\end{align*}
has a sequentially continuous differentiable solution $y=\phi(x)$ on the interval $-\alpha\leq x\leq \alpha$ where $\alpha=\min(a,b/M)$ and 
\[
M=\max\{|f(x,y)|\colon -a\leq x\leq a, \; -b\leq y\leq b\}
\]
\end{theorem}

In \cite[Section 2]{cheridito2015conditional}, some basic results in linear algebra are extended to the space $L^0(\mathbb{R})^n$ which culminates in a conditional version of the orthogonal decomposition theorem \cite[Corollary 2.12]{cheridito2015conditional}, whereby $L^0(\mathbb{R})^n$ is viewed as an module of rank $n$ over the commutative ring $L^0(\mathbb{R})$.  
By basic linear algebra in second-order arithmetic \cite{simpson2009subsystems}, $L^0(\mathbb{R})^n$ is a vector space of dimension $n$ in $\mathcal{S}$, and all results in \cite[Section 2]{cheridito2015conditional} are consequences of the transfer principle. 
Actually, they directly extend to the case where $n(\omega)$ is a measurable dimension.

\subsection{Metric spaces}\label{sec:metric} 

In the axiomatic system of second-order arithmetic, the only definable metric spaces are separable and complete ones which are coded as a completion of a countable set with a prescribed rate of convergence, see Simpson \cite[Definition II.5.1]{simpson2009subsystems}. 
We characterize a complete and separable metric space in the conditional model $\mathcal{S}$ as a vector metric space in a standard model as follows.
\begin{definition}\label{d:polishspace}
A non-empty set $H$ is said to be a \emph{conditional set}\footnote{See \cite[Definition 2.1]{drapeau2016algebra} for a formal definition.}, if there exists a restriction operation $|$ such that for every sequence $(x_k)$ in $H$ and every measurable partition $(A_k)$ there exists a unique element $x\in H$ such $x_k|A_k=x|A_k$ for all $k$. 
We name  this unique element a concatenation and denote it by $x=\sum_k x_k|A_k$. 

 Let $H$ be a conditional set. 
A function $d\colon H\times H\to L^0_{+}(\mathbb{R})$ is a \emph{conditional metric}, if   
\begin{itemize}
\item $d(\sum_k (x_k,y_k)|A_k)=\sum_{k} d(x_k,y_k)|A_k$ for all sequences $(x_k,y_k)$ in $H\times H$ and measurable partitions $(A_k)$,    
\item $d(x,y)=0$ if and only if $x=y$, 
\item $d(x,y)=d(y,x)$ for all $x,y\in H$, 
\item $d(x,z)\leq d(x,y)+d(y,z)$ for all $x,y,z \in H$. 
\end{itemize} 
A conditional metric space $H$ is said to be 
\begin{itemize}
\item \emph{conditionally separable}, if there exists a countable set $G\subset H$ such that for all $x\in H$ there is a sequence $(x_k)$ in $G$ and $n_1<n_2<\ldots$ in $L^0(\mathbb{N})$ such that $d(x_{n_m},x)\to 0$ almost everywhere where $x_{n}:=\sum_k x_k|\{n=k\}$ for $n\in L^0(\mathbb{N})$; 
\item \emph{conditionally complete}, if for every \emph{conditional Cauchy sequence}\footnote{That is, a sequence $(x_k)$ such that for all $r\in L^0_{++}(\Q)$ there exists $n\in L^0(\mathbb{N})$ with $d(x_m,x_p)<r$ for all $m,p\in L^0(\mathbb{N})$ with $m,p\geq n$.} $(x_k)$ there exists $x$ such that $d(x_k,x)\to 0$ almost everywhere.     
\end{itemize}
A conditional metric space is said to be \emph{conditionally Polish}, if it is conditionally separable and conditionally complete.  
\end{definition} 
\begin{examples}\label{exp:metricspaces}
\begin{itemize}
\item[a)] Let $(E,d)$ be a separable metric space with countable dense subset $S$.      
The metric $d$ extends to a vector metric on $L^0(E)$ with values in $L^0_+(\mathbb{R})$ via $d(x,y)(\omega):=d(x(\omega),y(\omega))$ almost everywhere. 
Let $G$ be the set of all constant functions with values in $S$.   
See Drapeau et al.~\cite[Sections 2 and 4]{DJK16} for a proof that $L^0(E)$ is conditionally complete.  
\item[b)] Suppose that $(\Omega,\mathcal{F},\mathbb{P})$ is a standard Borel probability space, and let $\mathcal{G}\subset \mathcal{F}$ be a sub-$\sigma$-algebra. For a non-negative random variable $x\colon \Omega\to \mathbb{R}$, define the extended conditional expectation $\mathbb{E}[ x | \mathcal{G}]:=\lim_{n\to \infty} \mathbb{E}[ x\wedge n| \mathcal{G}]$. 
For a real number $p\in[1,\infty)$, let $L^{p}(\mathcal{F}|\mathcal{G})$ denote the space of $\mathcal{F}$-measurable random variables such that $\mathbb{E}[|x|^p|\mathcal{G}]<\infty$ almost everywhere. The space $L^{p}(\mathcal{F}|\mathcal{G})$ can be represented in `tensorial' form as $L^0(\mathcal{G},\mathbb{R})\cdot L^p(\mathcal{F},\mathbb{R})$ within the space $L^0(\mathcal{F},\mathbb{R})$, see \cite[Proposition 1]{cerreia2016conditional}.  
The mapping $x\mapsto (\mathbb{E}[|x|^p|\mathcal{G}])^{1/p}$ defines a \emph{$\mathcal{G}$-conditional norm} on $L^{p}(\mathcal{F}|\mathcal{G})$, i.e.~a vector norm with values in $L^0_+(\mathcal{G},\mathbb{R})$ which equips $L^{p}(\mathcal{F}|\mathcal{G})$ with the structure of a topological $L^0(\mathcal{G},\mathbb{R})$-module, see the discussion in Filipovic et al.~\cite[Example 2.5]{filipovic2009separation} and Cerreia-Vioglio et al.~\cite[Section 8]{cerreia2016conditional}. 
The induced $\mathcal{G}$-stable vector metric is given by $d(x,y):=\mathbb{E}[|x-y|^p|\mathcal{G}]$. 
Let $S$ be a countable dense in $L^p(\mathcal{F},\mathbb{R})$, and let $G$ be the space of all products $x y$ where $x$ is a constant rational-valued random variable and $y\in S$.   
Conditional completeness follows from Vogelpoth \cite[Theorem 3.2.3]{vogelpoth2009phd}.   
Notice that here the underlying measure space is $(\Omega,\mathcal{G},\mathbb{P})$.  
\item[c)] Let $E$ be a standard Polish space, and let $F\colon \Omega\rightrightarrows E$ be an Effros measurable and closed-valued map. 
Then the set of almost everywhere selections of $F$ have the structure of a conditional Polish space.   
\end{itemize} 
\end{examples}
\begin{remark}
A conditional metric induces a standard topology on the \emph{set} $H$ given by a base $\{B(x,r)\colon x\in H, r\in L^0_+(\mathbb{Q})\}$.   
This standard topology on $H$ is in general not $\mathbb{R}$-metrizable.  
For several basic topological notions such as continuity and convergence, we have equivalences between their conditional and standard variants, see Drapeau et al.~\cite[Section 3]{drapeau2016algebra} for a systematic study.  
\end{remark}

We have the following classification result. 
\begin{proposition}\label{p:condmetric}
Every conditional Polish space corresponds uniquely to a Polish space in the conditional model $\mathcal{S}$. 
\end{proposition}
\begin{proof}
By definition, every Polish space in $\mathcal{S}$ is a conditional Polish space in the standard sense of Definition \ref{d:polishspace} (one can change from a conditionally countable family to a standard countable family by considering only the constant-valued indices). 
Conversely, let $H$ be a conditional Polish space, and let $G\subset H$ be countable and conditionally dense. 
Then each $x\in H$ can be identified with a conditional sequence $(x_n)_{n\in L^0(\N)}$ in the concatenation hull 
\[
\text{st}(G)=\left\{\sum_k x_k|A_k\colon (x_k) \text{ in } G, (A_k) \text{ measurable partition}\right\}
\]
such that $m<n$ implies $d(x_m,x_n)<2^{-m}$ for all $n,m\in L^0(\N)$.  
Indeed,  by a conditional version of the axiom of choice, see Drapeau et al.~\cite[Theorem 2.26]{drapeau2016algebra}, one finds a conditional sequence $(x_n)_{n\in L^0(\N)}$ such that $x_n\in O_n$ for all $n\in L^0(\N)$ where 
\[
O_n=\{y\in \text{st}(G)\colon d(x,y)<1/2^{n+1}\}
\]   
As $\text{st}(G)$ is countable from the perspective of $\mathcal{S}$, one concludes that $H$ is a Polish space in the model $ \mathcal{S}$. 
\end{proof}
We prove a characterization of closed sets in $\mathcal{S}$ which establishes a link to set-valued analysis.  
\begin{proposition}
Let $E$ be a standard Polish space, and let $W\subset L^0(E)$ be stable\footnote{Stability here and else where always refers to stability respectively closedness with respect to countable concatenations.}.  
Then the following are equivalent.  
\begin{itemize}
\item[(i)] $W$ is closed in $\mathcal{S}$.
\item[(ii)] $W$ is sequentially closed.   
\item[(iii)] There exists an Effros measurable and closed-valued map $X\colon \Omega\to 2^E$ such that $W$ coincides with the set of almost everywhere selections of $X$.   
\end{itemize}
\end{proposition}
\begin{proof}
We prove (i) $\Rightarrow$ (ii). 

Let $(x_k)$ be a sequence in $W$ such that $d(x_k,x)\to 0$ almost everywhere for some $x\in L^0(E)$, and let $G$ be a countable dense set in $E$. 
By contradiction, if $\mu(i(x,W^\sqsubset))>0$, then one finds a conditional ball $B(y,r)$ with center $y\in L^0(G)$ and radius $r\in L^0_{++}(\mathbb{Q})$  such that $\mu(i(x,B(y,r)\sqcap W^\sqsubset))>0$. 
But then there is also some $k$ such that $\mu(i(x_k,W^\sqsubset))>0$ which is the desired contradiction. 
Hence $\mu(i(x,W^\sqsubset))=0$, and therefore $x\in W$.   

We prove (ii) $\Rightarrow$ (i).  

Let 
\[
I:=\{(y,r)\in L^0(G)\times L^0_{++}(\mathbb{Q})\colon B(y,r)\cap W\neq \emptyset\}
\]
Notice that $I$ is conditionally countable.  
By a conditional version of the axiom of choice \cite[Theorem 2.26]{drapeau2016algebra}, one finds a conditional sequence $(x_i)$ such that $x_i\in B(y,r)\cap W$ for each $i=(y,r)\in I$.    
Now $W$ coincides with the Polish space defined by the countable set $\{x_i\}$ in the model $\mathcal{S}$.  

The equivalence of (ii) and (iii) is proved in Jamneshan and Zapata \cite[Section 5]{jamneshan2017compact}. 
\end{proof} 
As a direct consequence of the transfer principle, one obtains a conditional version of the \emph{Baire category theorem} (see Simpson, \cite[Theorem II.5.8]{simpson2009subsystems}), \emph{Urysohn's lemma} \cite[Lemma II.7.3]{simpson2009subsystems}, \emph{Tietze's extension theorem} \cite[Theorem II.7.5]{simpson2009subsystems}, the \emph{Ascoli lemma} \cite[Theorem III.2.8]{simpson2009subsystems},  the \emph{Heine-Borel theorem}  \cite[Theorem IV.1.5]{simpson2009subsystems}, a \emph{choice principle for compact sets} \cite[Theorem IV.1.8]{simpson2009subsystems}, and \emph{Tychonoff's theorem}  \cite[Theorem III.2.5]{simpson2009subsystems}.    
A conditional version of these theorems except Urysohn's lemma and Tietze's extension theorem were proved by hand previously, see Jamneshan and Zapata \cite{jamneshan2017compact} for references. 
We state a conditional version of Urysohn's lemma in standard language. 
Recall that $[0,1]$ denotes the set of all real-valued measurable functions $0\leq x\leq 1$. 
\begin{theorem}
Let $H$ be a conditional Polish space, and let $W_0$ and $W_1$ be stable and sequentially closed subsets of $H$ such that $W_0\sqcap W_1=\{\ast\}$.   
Then there exists a stable function $f\colon H\to [0,1]$ such that $f(x_k)
\to f(x)$ almost everywhere whenever $x_k\to x$ almost everywhere and $\{\omega\colon f(x)(\omega)=j\}=i(x,W_j)$ for all $j=0,1$ and $x\in H$. 
\end{theorem}
Heuristically, the previous theorem states that two random sets with values in a Polish space which only intersect on a negligible set can be separated by a random $\{0,1\}$-valued function $f$ in the sense that the largest measurable set on which a random element $x$ falls into the first/second set coincides almost surely with the measurable set $\{\omega\colon f(x)(\omega)=0/1\}$. 
A notable aspect about this 'randomized' Urysohn's lemma is that we do not leave the Borel measurable world without invoking a measurable selection argument, or in other words, measurability is a result of the construction. 
\begin{remark}
Conditional metric spaces were introduced in full generality in Drapeau et al.~\cite{drapeau2016algebra}. 
Conditional metric spaces were applied in Jamneshan et al.~\cite{jamneshan2017parameter} to stochastic control. 
It was argued there that they are a viable alternative to the established measurable selection techniques. 
Moreover, examples of problems were presented which can be solved with conditional analysis techniques which, however, are beyond the scope of applicability of measurable selection techniques due to topological restrictions.  
\end{remark} 

\subsection{Banach space theory}\label{sec:banach}

In this section, we develop basics of Banach spaces in $\mathcal{S}$, and connect it with functional analysis in $L^0(\mathbb{R})$-modules. 
An $L^0(\mathbb{R})$-module $H$ is a module over the commutative algebra $L^0(\mathbb{R})$. 
This yields a function from the underlying measure space $(\Omega,\mathcal{F},\mu)$ to $H$ by scalar multiplication with indicator functions $1_A\cdot x$, $A\in\mathcal{F}$ and $x\in H$\footnote{Actually, it is the measure algebra associated to $(\Omega,\mathcal{F},\mu)$ which acts on $H$.}. 
This enables to formalize concatenations in $H$:  
An element $x\in H$ is said to be a concatenation of a sequence $(x_k)$ in $H$ and a measurable partition $(A_k)$, if $1_{A_k}\cdot x = 1_{A_k} \cdot x_k$ for all $k$. 
We say that $H$ satisfies the countable concatenation property, if for every sequence $(x_k)$ in $H$ and each measurable partition $(A_k)$ there exists a unique concatenation in $H$. 
Throughout all $L^0(\mathbb{R})$-modules are assumed to satisfy the countable concatenation property.  
\begin{definition}
 Let $H$ be an $L^0(\mathbb{R})$-module. 
A function $\|\cdot\|\colon H\to L^0_+(\mathbb{R})$ is said to be a \emph{conditional norm}, if  
\begin{itemize}
\item $\|x\|=0$ if and only if $x=0$,
\item $\|r x\|=|r|\|x\|$ for all $r\in L^0(\mathbb{R})$ and $x\in H$, 
\item  $\|x+ y\|\leq \|x\|+\|y\|$ for all $x,y\in H$. 
\end{itemize}
\end{definition}
Since $H$ satisfies the countable concatenation property, a conditional norm is a stable function. 
Then $d(x,y):=\|x-y\|$ gives $H$ the structure of a conditional metric space. 
A conditional norm induces a standard topology given by the base 
$B(x,r):=\{y\in H\colon \|x-y\|<r\}$, $r\in L^0_{++}(\mathbb{Q})$ and $x\in H$.  
The generated topology renders addition and $L^0(\mathbb{R})$-scalar multiplication continuous where we consider on $L^0(\mathbb{R})$ its order topology, see remark \ref{r:euclidean}. 
From proposition \ref{p:condmetric} we have   
\begin{proposition}
There is a one-to-one correspondence from the class of conditionally separable and conditionally complete $L^0(\R)$-normed modules in a standard model and the class of separable Banach spaces\footnote{A separable Banach space in second-order arithmetic is introduced in \cite[Definition II.10.1]{simpson2009subsystems}.} 
 in the model $\mathcal{S}$. 
\end{proposition}    
One obtains the following examples of Banach spaces in the model $\mathcal{S}$ (see examples \ref{exp:metricspaces} for details). 
\begin{examples}
\begin{itemize}
\item[a)]  For $p\in [1,\infty)$, the conditional $L^p$-space $L^p(\mathcal{F}|\mathcal{G})$ with respect to a standard Borel probability space. 
\item[b)] The Bochner space $L^0(E)$ of equivalence classes of strongly measurable functions  with values in a separable normed vector space $E$.  
\end{itemize}
\end{examples}
As a consequence of the transfer principle, we obtain a conditional version of the 
 \emph{Hahn-Banach extension and separation theorems} (see Simpson \cite[Theorems IV.9.3, X.2.1]{simpson2009subsystems}), the \emph{Banach-Steinhaus theorem} \cite[Theorem II.10.8]{simpson2009subsystems}, the \emph{Banach-Alouglu theorem} \cite[Remark X.2.4]{simpson2009subsystems}, and the \emph{Krein-\v{S}mulian theorem} \cite[Theorem  X.2.7]{simpson2009subsystems} for conditionally separable Banach spaces. 
 We spell out a conditional version of the Banach-Steinhaus theorem below.  
The following definition is a standard interpretation of \cite[Definition II.10.5]{simpson2009subsystems} in $\mathcal{S}$. 
\begin{definition}
Let $H$ and $K$ be conditionally separable and conditionally complete $L^0(\R)$-normed modules. 
A function $f\colon H\to K$ is said to be 
\begin{itemize}
\item \emph{conditionally linear}, if $f(x + ry)=f(x)+ rf(y)$ for all $x,y\in H$ and $r\in L^0(\mathbb{R})$, 
\item \emph{conditionally bounded}, if there exists $r\in L^0_{++}(\mathbb{Q})$ such that $\|f(x)\|_K\leq r \|x\|_H$ for all $x\in H$.
\end{itemize}  
\end{definition}
Notice that any conditionally linear function is stable.  
One can deduce from \cite[Theorem II.10.7]{simpson2009subsystems} that a conditionally linear and conditionally bounded operator is sequentially continuous, i.e.~$\|f(x_k)-f(x)\|_K\to 0$ almost everywhere whenever $\|x_k-x\|_H\to 0$ almost everywhere. 
A conditional version of the Banach-Steinhaus theorem is obtained by an application of the transfer principle to \cite[Theorem II.10.8]{simpson2009subsystems}. 
One could think of its application in random operator theory and random differential equations, see e.g.~Skorohod~\cite{skorohod1983random} and Strand \cite{strand1970random}.

\begin{theorem}
Let $H$ and $K$ be conditionally separable and complete $L^0(\R)$-normed modules, and let $(f_n)$ be a standard sequence of $L^0(\mathbb{R})$-linear and sequentially continuous functions $f_n\colon H\to K$.  
If for every $x\in H$ there exists $r\in L^0_{++}(\mathbb{Q})$ such that $\|f_n(x)\|_K\leq r$ for all $n$, then there exists $q\in L^0_{++}(\mathbb{Q})$ such that $\|f_n(x)\|_K\leq q \|x\|_H$ for all $x$ and every $n$.  
\end{theorem}

\begin{remark}
Conditional extension and separation arguments are applied in Filipovic et al.~\cite{filipovic2012approaches} and Frittelli and Maggis \cite{frittelli2014complete} in risk measure theory.  
A conditional Fenchel-Moreau theorem is applied in Drapeau et al.~\cite{drapeau2017fenchel} in vector duality. 
A conditional Riesz representation theorem is used in Drapeau and Jamneshan \cite{drapeau2016conditional} in decision theory.  
\end{remark} 

\begin{remark}
Separation and duality results in topological $L^0(\R)$-modules are established in Guo et al.~\cite{guo2009separation} and Filipovic et al.~\cite{filipovic2009separation} with respect to two types of module topologies respectively, their connection is discussed in Guo \cite{guo2010relations}, see also the discussion in Jamneshan and Zapata \cite{jamneshan2017compact}.   
The standard topologies employed in \cite{guo2009separation}, and with a different motivation in Haydon et al.~\cite{haydon1991randomly}, are extensions of the topology of convergence in probability to $L^0(\mathbb{R})$-vector norms on general $L^0(\mathbb{R})$-modules. 
The class of such probabilistic topologies does not yield a comprehensive functional analytic discourse in $L^0(\R)$-modules, see Jamneshan and Zapata \cite{jamneshan2017compact}. 
  
The second type of topologies introduced in \cite{filipovic2009separation} proved to be more susceptible to a comprehensive functional analytic discourse for which strong evidence was provided in \cite{drapeau2016algebra} by embedding $L^0(\R)$-module theory in  conditional set theory and conditional topology.  
Conditional topological vector spaces are introduced in \cite[Section 5]{drapeau2016algebra}. 
In Drapeau et al.~\cite{DJK16} conditional completions of standard metric spaces are constructed and their connection to Lebesgue-Bochner spaces is established. 
The above listed consequences of the transfer principle have been proved by hand previously, see \cite{jamneshan2017compact} for an overview and references.  

A school of Russian mathematicians, starting with Kantorovic, studied to which extent results in functional analysis remain true if the real numbers are replaced by a Dedekind complete vector lattice of which $L^0(\mathbb{R})$ is one example. 
These investigations were naturally connected to Boolean-valued models. 
We refer to Kusraev and Kutateladze \cite{kusraev2012boolean} and Kutateladze \cite{kutateladze2012nonstandard} for an extensive overview of this tradition. 
\end{remark}


\subsection{Hilbert spaces} 

In this subsection, we will elaborate on basic results in separable Hilbert spaces\footnote{See \cite[Definition 9.3]{avigad2006fundamental} for a definition in second-order arithmetic.} where our focus lies on the conditional $L^2$-space 
\[
L^2(\mathcal{F}|\mathcal{G})=\{x\in L^0(\mathcal{F},\R)\colon \mathbb{E}[x^2|\mathcal{G}]<\infty\} 
\]
where $(\Omega,\mathcal{F},\mathbb{P})$ is a standard Borel probability space and $\mathcal{G}\subset \mathcal{F}$ is a sub-$\sigma$-algebra, see examples \ref{exp:metricspaces}. 
The inner product in $L^2(\mathcal{F}|\mathcal{G})$ is defined by $\langle x,y\rangle:= \mathbb{E}[xy|\mathcal{G}]$.  
Recall that for $L^2(\mathcal{F}|\mathcal{G})$ the base space is $(\Omega,\mathcal{G},\mathbb{P})$.

\begin{definition}
A sequence $(x_n)_{n\in L^0(\N)}$ in $L^2(\mathcal{F}|\mathcal{G})$  is said to be \emph{orthonormal}, if
\begin{itemize}
\item $\langle x_n,x_m\rangle=0$ whenever $\mathbb{P}(\{n=m\})=0$,
\item $\|x_n\|=1$ for all $n$.
\end{itemize}  
An orthonormal sequence $(x_n)_{n\in L^0(\mathbb{N})}$ is \emph{generating} if for every $x\in H$ there exists a sequence $(r_n)_{n\in L^0(\N)}$ in $L^0(\mathcal{G},\mathbb{R})$ such that $x=\sum_n r_n x_n$\footnote{Notice that this sum does not have a well-defined value in a standard setting a priori as it is uncountable, as a limit in $\mathcal{S}$ it is though meaningful, and this value can then be interpreted in a standard setting.}.  
A generating orthonormal sequence is called an \emph{orthonormal basis}.  
\end{definition}
The transfer principle applied to Avigad and Simic \cite[Theorem 10.9]{avigad2006fundamental} yields the following new result. 
\begin{theorem}\label{t:onb}
$L^2(\mathcal{F}|\mathcal{G})$ has an orthonormal basis in $\mathcal{S}$. 
\end{theorem}

\begin{remark}
The significance of theorem \ref{t:onb} lies in the fact that generally an $L^0(\R)$-module with the countable concatenation property has an algebraic basis if and only if it is finitely ranked, see Jamneshan and Zapata \cite[Proposition 3.5]{jamneshan2017compact}. 
This means that an infinite dimensional space such as $L^2(\mathcal{F}|\mathcal{G})$ does not have any module linear basis in a standard sense in general\footnote{For instance, if $(\Omega,\F,\mathbb{P})$ is purely non-atomic.}. 
The difference in $\mathcal{S}$ is that one allows the family which forms the basis to be a conditional family, cf.~\cite[Section 3]{jamneshan2017compact}. 
\end{remark}

We have the following \emph{projection theorem} thanks to the transfer principle and \cite[Theorem 12.5]{avigad2006fundamental}. 

\begin{theorem}
Let $W$ be a sequentially closed sub-module of $L^2(\mathcal{F}|\mathcal{G})$.   
Then every point $x\in L^2(\mathcal{F}|\mathcal{G})$ has a smallest distance $d(x,W):=\inf_{y\in W} d(x,y)$.   
\end{theorem}

A consequence of the projection theorem is the \emph{orthogonal decomposition theorem}: 

\begin{theorem}
Let $W$ be a sequentially closed sub-module of $L^2(\mathcal{F}|\mathcal{G})$.
Then there exists a sequentially closed sub-module $W^\bot$ such that each element $x\in L^2(\mathcal{F}|\mathcal{G})$ has a unique decomposition $x=z+y$ where $z\in W$ and $y\in W^\bot$. 
\end{theorem}
 
We have also a \emph{Riesz representation theorem} due to \cite[Theorem 13.4]{avigad2006fundamental} which interpreted in a standard setting reads as follows.  

\begin{theorem}
Let $f\colon L^2(\mathcal{F}|\mathcal{G}) \to L^0(\mathbb{R})$ be an $L^0(\R)$-linear and sequentially continuous function. 
Then there exists $y\in L^2(\mathcal{F}|\mathcal{G})$ such that $f(x)=\mathbb{E}[ xy |\mathcal{G}]$.   
\end{theorem} 
We obtain the following extension of \emph{von Neumann's mean ergodic theorem}, see \cite{avigad2006fundamental} for the result in ACA$_0$. 
\begin{theorem}
Let $T\colon L^2(\mathcal{F}|\mathcal{G})\to L^2(\mathcal{F}|\mathcal{G})$ be $L^0(\R)$-linear such that $\|Tx\|\leq \|x\|$ for all $x\in L^2(\mathcal{F}|\mathcal{G})$.   
Then $1/n(x + Tx +\ldots + T^{n-1} x)$ converges in the vector norm of $L^2(\mathcal{F}|\mathcal{G})$ almost everywhere to the projection of $x$ to the sequentially closed sub-module of all $T$-invariant vectors.  
\end{theorem} 

\begin{remark}
The conditional Hilbert space $L^2(\mathcal{F}|\mathcal{G})$ was introduced in Hansen and Richard \cite{hansen1987role} for purposes of financial modeling. 
Some applications in stochastic analysis of an orthogonal decomposition result are described in  Cerreia-Vioglio et al.~\cite{cerreia2017orthogonal}. 
A Riesz representation theorem, a projection theorem and an orthogonal decomposition theorem in complete random inner product spaces are proved in Guo \cite[Section 4]{guo2013homo}. 
A mean ergodic theorem for complete random inner product spaces is established in Guo and Zhang  \cite{guo2011neumann}. 
\end{remark}

\subsection{Fixed point theorems}

We close this section with fixed point theorems in the 'conditional model' $\mathcal{S}$. 
One could think of applications in equilibrium theory and random differential equations. 
 
A conditional version of the \emph{Brouwer fixed point theorem} was proved in Drapeau et al. \cite{drapeau2013brouwer}.  
The precise statement in the standard model is the following. 
\begin{theorem}
Let $x_1,\ldots,x_k$ be finitely many points in $L^0(\mathbb{R})^n$, $n\in \mathbb{N}$.  
Let 
\[
S=\Bigg\{\sum_{i=1}^k \lambda_i x_i\colon 0\leq \lambda_i\leq 1, i=1,\ldots,k\Bigg\}
\]
Let $f\colon S\to S$ be a stable sequentially continuous function.  
Then there is $x\in S$ such that $f(x)=x$.  
\end{theorem}
The proof in \cite{drapeau2013brouwer} is based on an adaptation of Sperner's lemma to a conditional setting.  
One obtains Sperner's lemma and the Brouwer fixed point theorem in the following slightly more general form as a consequence of the transfer principle applied to Simpson \cite[Theorem IV.7.6]{simpson2009subsystems}.
\begin{theorem}
Let $k,n\in L^0(\N)$ and $(x_i)_{i\leq k}$ be a conditionally finite sequence of points in $L^0(\mathbb{R})^n$. 
Without loss of generality suppose that $k=\sum_m k_m|A_m$ and $n=\sum_m n_m|A_m$, and let $(x^m_i)_{i\leq k_m}$ be finitely many elements in $L^0(\mathbb{R})^{n_m}$ which comprise the conditional sequence $(x_i)_{i\leq k}$ on $A_m$ for each $m$. 
Build the conditional simplex
\[
S:=\sum_m \Bigg\{\sum_{i=1}^{k_m} \lambda_i x^m_i\colon 0\leq \lambda_i\leq 1, i=1,\ldots,k_m\Bigg\}\Bigg|A_m
\]  
Then every stable sequentially continuous function $f\colon S\to S$ has a fixed point.  
\end{theorem}
As a novel result, we deduce a conditional version of the \emph{Markov-Kakutani fixed point theorem} by applying the transfer principle to \cite[Lemma IV.9.1]{simpson2009subsystems}. 
\begin{theorem}
Let  
\[
[-1,1]^{L^0(\N)}:=\{(x_k)_{k\in L^0(\N)}\colon x_k\in L^0(\mathbb{R}),\; -1\leq x_k\leq 1 \text{ for all } k\}
\]
endowed with the conditional metric\footnote{Recall that the sum is uncountable in a standard world, but makes sense in a 'conditional model' as the limit of conditionally finite sums.}  
\[
d((x_k),(y_k)):=\sum_k 2^{-k}\frac{|x_k-y_k|}{|x_k-y_k|+1}
\]
Let $S$ be a stable, sequentially closed and $L^0(\mathbb{R})$-convex\footnote{This means $\lambda x + (1-\lambda)y\in S$ whenever $x,y\in S$ and $\lambda\in L^0(\mathbb{R})$ with $0\leq \lambda \leq 1$.} subset of $[-1,1]^{L^0(\N)}$.  
Let $(f_k)$ be a conditional sequence\footnote{This means for functions that $\big(\sum_n f_{k_n}|A_n\big)(x)=\sum_n f_{k_n}(x)|A_n$ for all $x$, every sequence $(k_n)$ in $L^0(\mathbb{N})$ and measurable partition $(A_n)$.} of $L^0(\mathbb{R})$-affine and sequentially continuous functions $f_k\colon S\to S$ such that $f_k\circ f_m=f_m\circ f_k$ for all $m,k\in L^0(\N)$.  
Then there exists $x\in S$ such that $f_k(x)=x$ for all $k$. 
\end{theorem}


\section{Further connections and outlook}\label{s5}
We have left out measure theory.  
For the development of measure theory in second-order arithmetic, we refer to Simic \cite{simic2004thesis}, Simpson \cite{simpson2009subsystems} and Yu \cite{yu1993riesz}.  
In conditional set theory, basic results in measure theory are established in Jamneshan et al.~\cite{jamneshan2017measures}.  
As measure theory in second-order arithmetic is based on Daniell's functional approach to integration, the conditional version of the Daniell-Stone theorem and 
Riesz representation theorem in \cite[Section 5]{jamneshan2017measures} can be connected to respective theorems in $\mathcal{S}$. 
In \cite[Section 4]{jamneshan2017measures}, a connection between kernels in a standard model and measures in conditional set theory are proved.  
This connection fully carries over to the model $\mathcal{S}$, for Borel probability measures on $H$ where $H$ is a conditionally compact metric space. 
Perspectively, understanding the standard meaning of the existence of a Haar measure (\cite[Section X.1]{simpson2009subsystems}), 
and of the maximal ergodic and pointwise ergodic theorems (\cite[Chapter 5]{simic2004thesis}) in $\mathcal{S}$ may be of interest.  

For some mathematical theorems relevant for applications in mathematical economics, stronger second-order axioms than ACA$_0$ are required. 
For example, weak$^\ast$-closures play an important role in a recent duality result in vector optimization in Grad and Jamneshan \cite{grad2018perturbational} which is useful for applications such as minimization of conditional risk measures. Now, by \cite[Theorem X.2.9]{simpson2009subsystems}, $\Pi_{1}^{1}$-CA$_{0}$ is required to prove that weak$^{\ast}$-closures in the normal dual of 
a separable Banach space even exist. Moreover, it has been argued e.g.~by Kohlenbach in \cite{kohlenbach2000higher} that the second-order setting, comprehensive as it is, 
is not sufficient for important areas like functional analysis or topology, and should therefore be extended to higher-order frameworks that also allow for quantification over sets 
of real numbers, sets of sets of real numbers, sets thereof etc.

It would thus be worthwhile (1) to consider the validity of stronger comprehension axioms like $\Pi_{1}^{1}$-CA$_{0}$ (see \cite{simpson2009subsystems}) or, more generally, $\Pi_{1}^{n}$-CA$_{0}$ 
in our conditional model, and (2) to come up with conditional interpretations of higher-order systems and a corresponding transfer principle.
The question whether $\Pi_{1}^{1}$-CA$_{0}$ holds in our model is not immediately answered by the method used above for ACA$_{0}$. For this purpose, it would be necessary to consider projections of sets of real numbers, which would require working in a conditional power set of the conditional real numbers. 

Nevertheless, we are optimistic that such extensions can be constructed. Indeed, following the general construction of a conditional 	power set in \cite{drapeau2016algebra}, more precisely its variant defined in \cite{jamneshan2017measures}, the conditional power set of the conditional power set of the conditional natural 
numbers can be formed, which indicates the necessary construction needed to prove the validity of stronger comprehension schemes and  a transfer principle for higher-order systems.
We plan to develop this in detail in future work.

\end{document}